\documentclass{article}
\usepackage{amssymb}
\usepackage{amsfonts}
\usepackage{amsmath}

\input{tcilatex}
\begin{document}

\title{Compact embeddings and Pitt's property for weighted sequence spaces
of Sobolev type}
\author{Pierre-A. Vuillermot$^{\ast ,\ast \ast }$ \\
%EndAName
Center for Mathematical Studies, CEMS.UL, Faculdade de Ci\^{e}ncias,\\
Universidade de Lisboa, 1749-016 Lisboa, Portugal$^{\ast }$\\
Universit\'{e} de Lorraine, CNRS, IECL, F-54000 Nancy, France$^{\ast \ast }$%
\\
pavuillermot@fc.ul.pt\\
ORCID Number: 0000-0003-3527-9163}
\date{}
\maketitle

\begin{abstract}
In this article we introduce a new class of weighted sequence spaces of
Sobolev type and prove several compact embedding theorems for them. It is
our contention that the chosen class is general enough so as to allow
applications in various areas of mathematics and mathematical physics. In
particular, our results constitute a generalization of those compact
embeddings recently obtained in relation to the spectral analysis of a class
of master equations with non-constant coefficients arising in
non-equilibrium statistical mechanics. As a byproduct of our considerations,
we also prove a theorem of Pitt's type asserting that under some conditions
every linear bounded transformation from one weighted sequence space of the
class into another is compact.



\ \ \ \ \ \ \ \ \ \ \ \textbf{Abbreviated title}: Compact Embeddings
\end{abstract}

\section{Introduction and outline}

The essential role played by compact embeddings of Sobolev spaces of various
kinds in the analysis of initial- and boundary-value problems for ordinary
and partial differential equations is well known (see, e.g., \cite%
{adamsfournier} or \cite{edmundsevans} and the numerous references therein).
Of equal importance are certain Hilbert spaces of sequences and their
relation to Sobolev spaces of periodic functions as in \cite{bisgard}. In
Section 2 of this article we introduce a new scale of weighted sequence
spaces of Sobolev type and prove compact embedding results for them. The
chosen class is general enough so as to allow applications in various areas
of mathematics or mathematical physics. In particular, our results
generalize those embedding properties recently used in \cite{vuillermot} in
relation to the analysis of a class of master equations with non-constant
coefficients arising in non-equilibrium statistical mechanics, thereby
extending the investigations started in \cite{boglivuillermot0} and \cite%
{boglivuillermot1}. We also prove there a theorem of Pitt's type asserting
that under some restrictions, every linear bounded transformation between
two spaces of the scale is compact. We refer the reader to \cite{pitt} for
the original statement involving a linear bounded operator between two
spaces of summable sequences, to \cite{delpech} and \cite{fabianzizler} for
much shorter proofs thereof and to Theorem 2.1.4 in \cite{albiackalton} or
Proposition 2.c.3 in \cite{lindenstrausstza} for yet more condensed
arguments.

\section{The results}

With $s\in \left[ 1,+\infty \right) $, $k\in \mathbb{R}$ and $w=(w_{\mathsf{m%
}})_{\mathsf{m}\in \mathbb{Z}}$ a sequence of weights satisfying $w_{\mathsf{%
m}}>0$ for every $\mathsf{m}$, let us consider the separable Banach space $%
h_{\mathbb{C},w}^{k,s}$ of Sobolev type consisting of all complex sequences $%
\mathsf{p}=(p_{\mathsf{m}})$ endowed with the usual algebraic operations and
the norm%
\begin{equation}
\left\Vert \mathsf{p}\right\Vert _{k,s,w}:=\left( \dsum\limits_{\mathsf{m\in 
}\mathbb{Z}}w_{\mathsf{m}}\left( 1+\left\vert \mathsf{m}\right\vert
^{s}\right) ^{k}\left\vert p_{\mathsf{m}}\right\vert ^{s}\right) ^{\frac{1}{s%
}}<+\infty .  \label{norm}
\end{equation}%
If $k=0$ we simply write $l_{\mathbb{C},w}^{s}:=h_{\mathbb{C},w}^{0,s}$ and 
\begin{equation}
\left\Vert \mathsf{p}\right\Vert _{s,w}:=\left( \dsum\limits_{\mathsf{m\in }%
\mathbb{Z}}w_{\mathsf{m}}\left\vert p_{\mathsf{m}}\right\vert ^{s}\right) ^{%
\frac{1}{s}}  \label{normbis}
\end{equation}%
for the corresponding norm. We may refer to $s$ as the \textit{degree of
summability }of $\mathsf{p}$ and to $k$ as its \textit{generalized order of
differentiability}, a terminology justified by the analogy with the usual
Sobolev space theory and its relation to Fourier analysis on $\mathbb{R}^{d}$
(see, e.g., Chapter VI in \cite{yosida}). Unless $s=1$, the spaces $h_{%
\mathbb{C},w}^{k,s}$ are reflexive and it is also easily determined that $h_{%
\mathbb{C},w}^{k,s}$ is a Hilbert space if and only if $s=2$, in which case (%
\ref{norm}) is related to the sesquilinear form%
\begin{equation*}
\left( \mathsf{p,q}\right) _{k,2,w}:=\dsum\limits_{\mathsf{m\in }\mathbb{Z}%
}w_{\mathsf{m}}\left( 1+\mathsf{m}^{2}\right) ^{k}p_{\mathsf{m}}\bar{q}_{%
\mathsf{m}}
\end{equation*}%
in the usual way. We also note that the role of the sequence of weights $(w_{%
\mathsf{m}})$ in the Hilbert space case generally speaking amounts to making
certain initial-value problems self-adjoint (see our discussion following
the proof of Theorem 2 below).

Our proofs of various compactness criteria for the space $h_{\mathbb{C}%
,w}^{k,s}$ rest in an essential way on the existence of a very concrete and
simple Schauder basis therein (for the definition and many properties of
such bases see, e.g., \cite{albiackalton} and \cite{lindenstrausstza}).
Thus, for any $\mathsf{m}\in \mathbb{Z}$ we define $\mathsf{e}_{\mathsf{m}}$
by $\left( \mathsf{e}_{\mathsf{m}}\right) _{\mathsf{n}}=\delta _{\mathsf{m,n}%
}$ for every $\mathsf{n}\in \mathbb{Z}$ so that $\mathsf{e}_{\mathsf{m}}\in $
$h_{\mathbb{C},w}^{k,s}$ for each $\mathsf{m}$ with 
\begin{equation*}
\left\Vert \mathsf{e}_{\mathsf{m}}\right\Vert _{_{k,s,w}}=w_{\mathsf{m}}^{%
\frac{1}{s}}\left( 1+\left\vert \mathsf{m}\right\vert ^{s}\right) ^{\frac{k}{%
s}}.
\end{equation*}%
Cauchy's criterion then shows that every $\mathsf{p=}(p_{\mathsf{m}})\in h_{%
\mathbb{C},w}^{k,s}$ may be expanded in a unique way as the norm-convergent
series%
\begin{equation}
\mathsf{p=}\dsum\limits_{\mathsf{m}\in \mathbb{Z}}p_{\mathsf{m}}\mathsf{e}_{%
\mathsf{m}},  \label{expansion}
\end{equation}%
a fact that will be crucial in what follows. As a preliminary remark we note
that if $k\leq k^{\prime }$, there exists the continuous embedding%
\begin{equation}
h_{\mathbb{C},w}^{k^{\prime },s}\rightarrow h_{\mathbb{C},w}^{k,s}
\label{contembedding}
\end{equation}%
by virtue of the obvious inequality 
\begin{equation*}
\left\Vert \mathsf{p}\right\Vert _{k,s,w}^{s}=\dsum\limits_{\mathsf{m\in }%
\mathbb{Z}}w_{\mathsf{m}}\left( 1+\left\vert \mathsf{m}\right\vert
^{s}\right) ^{k}\left\vert p_{\mathsf{m}}\right\vert ^{s}\leq \dsum\limits_{%
\mathsf{m\in }\mathbb{Z}}w_{\mathsf{m}}\left( 1+\left\vert \mathsf{m}%
\right\vert ^{s}\right) ^{k^{\prime }}\left\vert p_{\mathsf{m}}\right\vert
^{s}=\left\Vert \mathsf{p}\right\Vert _{k^{\prime },s,w}^{s}.
\end{equation*}

Our first result is then:

\bigskip

\textbf{Theorem 1. }\textit{The following statements hold: }

\textit{(a) If the strict inequality }$k<k^{\prime }$\textit{\ holds, then
embedding (\ref{contembedding}) is compact, in which case we write}%
\begin{equation}
h_{\mathbb{C},w}^{k^{\prime },s}\hookrightarrow h_{\mathbb{C},w}^{k,s}.
\label{compactembedding}
\end{equation}

\textit{(b) Let us assume that }$0\leq k<k^{\prime }$\textit{, and that the
sequence of weights }$w=(w_{\mathsf{m}})$ \textit{satisfies the additional
constraint}%
\begin{equation}
\inf_{\mathsf{m\in }\mathbb{Z}}w_{\mathsf{m}}>0\text{.\label{infimum}}
\end{equation}%
\textit{Then if }$s\geq t\geq 1$, \textit{there exists the compact embedding}%
\begin{equation}
h_{\mathbb{C},w}^{k^{\prime },t}\hookrightarrow h_{\mathbb{C},w}^{k,s}.
\label{embeddingter}
\end{equation}%
\textit{In particular, the embedding}%
\begin{equation*}
h_{\mathbb{C},w}^{k^{\prime },t}\hookrightarrow l_{\mathbb{C},w}^{s}
\end{equation*}%
\textit{is compact.}

\bigskip

\textbf{Proof. }Let $\mathcal{K}$ be a bounded set in $h_{\mathbb{C}%
,w}^{k^{\prime },s}$ and let $\kappa >0$ be the radius of a ball centered at
the origin of $h_{\mathbb{C},w}^{k^{\prime },s}$ and containing $\mathcal{K}$%
. Then for each $\varepsilon >0$ there exists $\mathsf{m}^{\ast }\in \mathbb{%
Z}^{+}$ such that for every $\left\vert \mathsf{m}\right\vert \mathsf{\geq m}%
^{\ast }$ we have%
\begin{equation*}
\frac{1}{\left( 1+\left\vert \mathsf{m}\right\vert ^{s}\right) ^{k^{\prime
}-k}}\leq \left( \frac{\varepsilon }{2\kappa }\right) ^{s}
\end{equation*}%
since $k^{\prime }-k>0$. Therefore, for every $\mathsf{p}\in \mathcal{K}$ we
obtain%
\begin{eqnarray*}
&&\dsum\limits_{\left\vert \mathsf{m}\right\vert \mathsf{\geq m}^{\ast }}w_{%
\mathsf{m}}\left( 1+\left\vert \mathsf{m}\right\vert ^{s}\right)
^{k}\left\vert p_{\mathsf{m}}\right\vert ^{s} \\
&=&\dsum\limits_{\left\vert \mathsf{m}\right\vert \mathsf{\geq m}^{\ast }}%
\frac{1}{\left( 1+\left\vert \mathsf{m}\right\vert ^{s}\right) ^{k^{\prime
}-k}}w_{\mathsf{m}}\left( 1+\left\vert \mathsf{m}\right\vert ^{s}\right)
^{k^{\prime }}\left\vert p_{\mathsf{m}}\right\vert ^{s} \\
&\leq &\left( \frac{\varepsilon }{2\kappa }\right) ^{s}\dsum\limits_{\mathsf{%
m\in }\mathbb{Z}}w_{\mathsf{m}}\left( 1+\left\vert \mathsf{m}\right\vert
^{s}\right) ^{k^{\prime }}\left\vert p_{\mathsf{m}}\right\vert ^{s}\leq
\left( \frac{\varepsilon }{2\kappa }\right) ^{s}\kappa ^{s}=\left( \frac{%
\varepsilon }{2}\right) ^{s},
\end{eqnarray*}%
that is, 
\begin{equation}
\left( \dsum\limits_{\left\vert \mathsf{m}\right\vert \mathsf{\geq m}^{\ast
}}w_{\mathsf{m}}\left( 1+\left\vert \mathsf{m}\right\vert ^{s}\right)
^{k}\left\vert p_{\mathsf{m}}\right\vert ^{s}\right) ^{\frac{1}{s}}\leq 
\frac{\varepsilon }{2}.  \label{estimate}
\end{equation}%
Now let $\mathcal{\bar{K}}$ be the closure of $\mathcal{K}$ in $h_{\mathbb{C}%
,w}^{k,s}$. In order to show that $\mathcal{\bar{K}}$ is compact it is
necessary and sufficient to show that $\mathcal{\bar{K}}$ is totally bounded
or, equivalently, that there exists a finite-dimensional subspace $\mathcal{S%
}_{\varepsilon ,\kappa }\subset h_{\mathbb{C},w}^{k,s}$ such that the
distance of every $\mathsf{p}\in \mathcal{\bar{K}}$ to $\mathcal{S}%
_{\varepsilon ,\kappa }$ satisfies 
\begin{equation*}
\func{dist}\left( \mathsf{p,}\mathcal{S}_{\varepsilon ,\kappa }\right)
:=\inf_{\mathsf{q}\in \mathcal{S}_{\varepsilon ,\kappa }}\left\Vert \mathsf{%
p-q}\right\Vert _{k,s,w}\leq \varepsilon
\end{equation*}%
(see, e.g., Proposition 2.1 in \cite{bisgard}, or \cite{kirillovgvishiani}).
To this end we choose 
\begin{equation*}
\mathcal{S}_{\varepsilon ,\kappa }=\limfunc{span}\left\{ \mathsf{e}_{-%
\mathsf{m}^{\ast }},...,\mathsf{e}_{\mathsf{m}^{\ast }}\right\}
\end{equation*}%
where $\mathsf{e}_{\mathsf{m}}$ is as in (\ref{expansion}). Furthermore,
having (\ref{expansion}) in mind we define%
\begin{equation*}
\mathsf{p}_{\mathsf{m}^{\ast }}:\mathsf{=}\dsum\limits_{\left\vert \mathsf{m}%
\right\vert <\mathsf{m}^{\ast }}p_{\mathsf{m}}\mathsf{e}_{\mathsf{m}}.
\end{equation*}%
Then $\mathsf{p}_{\mathsf{m}^{\ast }}\in \mathcal{S}_{\varepsilon ,\kappa }$
and for every $\mathsf{p}\in $ $\mathcal{K}$ we obtain%
\begin{equation}
\func{dist}\left( \mathsf{p,}\mathcal{S}_{\varepsilon ,\kappa }\right) \leq
\left\Vert \mathsf{p-p}_{\mathsf{m}^{\ast }}\right\Vert _{k,s,w}=\left(
\dsum\limits_{\left\vert \mathsf{m}\right\vert \mathsf{\geq m}^{\ast }}w_{%
\mathsf{m}}\left( 1+\left\vert \mathsf{m}\right\vert ^{s}\right)
^{k}\left\vert p_{\mathsf{m}}\right\vert ^{s}\right) ^{\frac{1}{s}}\leq 
\frac{\varepsilon }{2}  \label{estimatequarto}
\end{equation}%
according to (\ref{estimate}). Finally, for every $\mathsf{p}\in \mathcal{%
\bar{K}}\diagdown \mathcal{K}$ there exists $\mathsf{p}_{\varepsilon }\in 
\mathcal{K}$ such that $\left\Vert \mathsf{p-p}_{\varepsilon }\right\Vert
_{k,s,w}\leq \frac{\varepsilon }{2}$ so that we obtain%
\begin{equation*}
\func{dist}\left( \mathsf{p,}\mathcal{S}_{\varepsilon ,\kappa }\right) \leq 
\frac{\varepsilon }{2}+\inf_{\mathsf{q}\in \mathcal{S}_{\varepsilon ,\kappa
}}\left\Vert \mathsf{p}_{\varepsilon }\mathsf{-q}\right\Vert _{k,s,w}\leq
\varepsilon
\end{equation*}%
as a consequence of (\ref{estimatequarto}), which proves Statement (a).

In order to prove Statement (b), we only consider $s>t$ since the case $s=t$
has already been dealt with. Since $\frac{s}{t}>1$ we first have%
\begin{eqnarray}
&&\dsum\limits_{\mathsf{m\in }\mathbb{Z}}w_{\mathsf{m}}\left( 1+\left\vert 
\mathsf{m}\right\vert ^{s}\right) ^{k}\left\vert p_{\mathsf{m}}\right\vert
^{s}  \notag \\
&=&\dsum\limits_{\mathsf{m\in }\mathbb{Z}}\left( w_{\mathsf{m}}^{^{\frac{t}{s%
}}}\left( 1+\left( \left\vert \mathsf{m}\right\vert ^{t}\right) ^{\frac{s}{t}%
}\right) ^{\frac{kt}{s}}\left\vert p_{\mathsf{m}}\right\vert ^{t}\right) ^{%
\frac{s}{t}}  \label{inequality} \\
&\leq &\left( \dsum\limits_{\mathsf{m\in }\mathbb{Z}}w_{\mathsf{m}}^{^{\frac{%
t}{s}}}\left( 1+\left( \left\vert \mathsf{m}\right\vert ^{t}\right) ^{\frac{s%
}{t}}\right) ^{\frac{kt}{s}}\left\vert p_{\mathsf{m}}\right\vert ^{t}\right)
^{\frac{s}{t}}\leq \left( \dsum\limits_{\mathsf{m\in }\mathbb{Z}}w_{\mathsf{m%
}}^{^{\frac{t}{s}}}\left( 1+\left\vert \mathsf{m}\right\vert ^{t}\right)
^{k}\left\vert p_{\mathsf{m}}\right\vert ^{t}\right) ^{\frac{s}{t}}  \notag
\end{eqnarray}%
since $k\geq 0$. Furthermore, from (\ref{infimum}) and the fact that $\frac{t%
}{s}<1$ we get%
\begin{equation}
w_{\mathsf{m}}^{^{\frac{t}{s}}}=\left( \inf_{\mathsf{m\in }\mathbb{Z}}w_{%
\mathsf{m}}\right) ^{\frac{t}{s}}\left( \frac{w_{\mathsf{m}}}{\inf_{\mathsf{%
m\in }\mathbb{Z}}w_{\mathsf{m}}}\right) ^{\frac{t}{s}}\leq \left( \inf_{%
\mathsf{m\in }\mathbb{Z}}w_{\mathsf{m}}\right) ^{\frac{t}{s}-1}w_{\mathsf{m}%
}.  \label{inequalitybis}
\end{equation}%
The substitution of (\ref{inequalitybis}) into the very last term of (\ref%
{inequality}) then gives%
\begin{equation*}
\dsum\limits_{\mathsf{m\in }\mathbb{Z}}w_{\mathsf{m}}\left( 1+\left\vert 
\mathsf{m}\right\vert ^{s}\right) ^{k}\left\vert p_{\mathsf{m}}\right\vert
^{s}\leq \left( \inf_{\mathsf{m\in }\mathbb{Z}}w_{\mathsf{m}}\right) ^{1-%
\frac{s}{t}}\left( \dsum\limits_{\mathsf{m\in }\mathbb{Z}}w_{\mathsf{m}%
}\left( 1+\left\vert \mathsf{m}\right\vert ^{t}\right) ^{k}\left\vert p_{%
\mathsf{m}}\right\vert ^{t}\right) ^{\frac{s}{t}},
\end{equation*}%
that is%
\begin{equation*}
\left\Vert \mathsf{p}\right\Vert _{k,s,w}\leq c_{s,t}\left\Vert \mathsf{p}%
\right\Vert _{k,t,w},
\end{equation*}%
for some finite constant $c_{s,t}>0$, so that the embedding%
\begin{equation}
h_{\mathbb{C},w}^{k,t}\rightarrow h_{\mathbb{C},w}^{k,s}
\label{continuousembedingbis}
\end{equation}%
is continuous. Therefore, (\ref{embeddingter}) may be viewed as the
composition%
\begin{equation*}
h_{\mathbb{C},w}^{k^{\prime },t}\hookrightarrow h_{\mathbb{C}%
,w}^{k,t}\rightarrow h_{\mathbb{C},w}^{k,s}
\end{equation*}%
where the first embedding is compact according to Statement (a). \ \ $%
\blacksquare $

\bigskip

It is interesting to note the interlacing properties of $l_{\mathbb{C}%
,w}^{s} $ with respect to $h_{\mathbb{C},w}^{k^{\prime },s}$ and $h_{\mathbb{%
C},w}^{k,s}$ as functions of the signs of $k$ and $k^{\prime }$:

\bigskip

\textbf{Corollary 1.}\textit{\ The following statements hold:}

\textit{(a) If }$k<k^{\prime }\leq 0$\textit{, we have the embeddings}%
\begin{equation}
l_{\mathbb{C},w}^{s}\rightarrow h_{\mathbb{C},w}^{k^{\prime
},s}\hookrightarrow h_{\mathbb{C},w}^{k,s}  \label{embedding}
\end{equation}%
\textit{where the first one is continuous and the second one compact. In
particular, the embedding}%
\begin{equation*}
l_{\mathbb{C},w}^{s}\hookrightarrow h_{\mathbb{C},w}^{k,s}
\end{equation*}%
\textit{is compact.}

\textit{(b) If }$k\leq 0<k^{\prime }$, \textit{the above chain of embeddings
is changed into}%
\begin{equation}
h_{\mathbb{C},w}^{k^{\prime },s}\hookrightarrow l_{\mathbb{C}%
,w}^{s}\rightarrow h_{\mathbb{C},w}^{k,s}  \label{embeddingbis}
\end{equation}%
\textit{where the first embedding is compact and the second one continuous.}

\textit{(c) If }$0\leq k<k^{\prime }$,\textit{\ the chain of embeddings
becomes}%
\begin{equation*}
h_{\mathbb{C},w}^{k^{\prime },s}\hookrightarrow h_{\mathbb{C}%
,w}^{k,s}\rightarrow l_{\mathbb{C},w}^{s}
\end{equation*}%
\textit{where the first embedding is compact and the second one continuous.}

\bigskip

\textbf{Proof. }The continuity of the first embedding in (\ref{embedding})
is a consequence of the inequality%
\begin{equation*}
\left\Vert \mathsf{p}\right\Vert _{k^{\prime },s,w}^{s}=\dsum\limits_{%
\mathsf{m\in }\mathbb{Z}}\frac{w_{\mathsf{m}}}{\left( 1+\left\vert \mathsf{m}%
\right\vert ^{s}\right) ^{\left\vert k^{\prime }\right\vert }}\left\vert p_{%
\mathsf{m}}\right\vert ^{s}\leq \dsum\limits_{\mathsf{m\in }\mathbb{Z}}w_{%
\mathsf{m}}\left\vert p_{\mathsf{m}}\right\vert ^{s}=\left\Vert \mathsf{p}%
\right\Vert _{s,w}^{s},
\end{equation*}%
so that the remaining part of Statement (a) then follows from Theorem 1.

As for Statement (b), the compactness of the first embedding in (\ref%
{embeddingbis}) follows from (\ref{compactembedding}) with $k=0$ while the
remaining part of the statement is a consequence of what has already been
proved, as is the case for Statement (c). \ \ $\blacksquare $

\bigskip

Our next theorem is concerned with embedding properties between spaces whose
norms involve different weights. A particular case of that result played a
central role in \cite{vuillermot} in relation to the analysis of a class of
master equations with non-constant coefficients arising in non-equilibrium
statistical mechanics. Thus, aside from $w$ let us introduce another
sequence of positive weights $\hat{w}=(\hat{w}_{\mathsf{m}})$ related to $w$
by the inequalities 
\begin{equation}
c_{1}\leq \frac{w_{\mathsf{m}}^{\frac{t}{s}}}{\hat{w}_{\mathsf{m}}}\leq c_{2}
\label{weightsbis}
\end{equation}%
for every $\mathsf{m}\in \mathbb{Z}$, where the constants $c_{1,2}>0$ are
finite, independent of $\mathsf{m}$ and $s,t\in \left[ 1,+\infty \right) $
as before. Then we have:

\bigskip

\textbf{Theorem 2.} \textit{Let us assume that (\ref{weightsbis}) holds with 
}$s>t$,\textit{\ and that }$k\in \left( \frac{s-t}{st},+\infty \right) $. 
\textit{Then\ there exist the} \textit{embeddings}%
\begin{equation}
h_{\mathbb{C},w}^{k,s}\hookrightarrow l_{\mathbb{C},\hat{w}}^{t}\rightarrow
l_{\mathbb{C},w}^{s}  \label{embeddings}
\end{equation}%
\textit{where the first one is compact and the second one continuous.}

\bigskip

\textbf{Proof.} We begin by proving the continuity of the first embedding in
(\ref{embeddings}). Let us set%
\begin{equation}
r:=\frac{st}{s-t}  \label{duality}
\end{equation}%
so that 
\begin{equation*}
\frac{1}{r}+\frac{1}{s}=\frac{1}{t}
\end{equation*}%
and $r\in \left( 1,+\infty \right) $. Owing to H\"{o}lder's inequality and
the first inequality in (\ref{weightsbis}) we then get%
\begin{eqnarray}
&&\dsum\limits_{\mathsf{m}\in \mathbb{Z}}\hat{w}_{\mathsf{m}}\left\vert p_{%
\mathsf{m}}\right\vert ^{t}  \notag \\
&=&\dsum\limits_{\mathsf{m}\in \mathbb{Z}}\left( \frac{1}{\left(
1+\left\vert \mathsf{m}\right\vert ^{s}\right) ^{\frac{k}{s}}}\right)
^{t}\left( \hat{w}_{\mathsf{m}}^{\frac{1}{t}}\left( 1+\left\vert \mathsf{m}%
\right\vert ^{s}\right) ^{\frac{k}{s}}\left\vert p_{\mathsf{m}}\right\vert
\right) ^{t}  \label{estimates} \\
&\leq &\left( \dsum\limits_{\mathsf{m}\in \mathbb{Z}}\frac{1}{\left(
1+\left\vert \mathsf{m}\right\vert ^{s}\right) ^{\frac{kr}{s}}}\right) ^{%
\frac{t}{r}}\left( \dsum\limits_{\mathsf{m}\in \mathbb{Z}}\hat{w}_{\mathsf{m}%
}^{\frac{s}{t}}\left( 1+\left\vert \mathsf{m}\right\vert ^{s}\right)
^{k}\left\vert p_{\mathsf{m}}\right\vert ^{s}\right) ^{\frac{t}{s}}  \notag
\\
&\leq &c_{1}^{-1}\left( \dsum\limits_{\mathsf{m}\in \mathbb{Z}}\frac{1}{%
\left( 1+\left\vert \mathsf{m}\right\vert ^{s}\right) ^{\frac{kr}{s}}}%
\right) ^{\frac{t}{r}}\left( \dsum\limits_{\mathsf{m}\in \mathbb{Z}}w_{%
\mathsf{m}}\left( 1+\left\vert \mathsf{m}\right\vert ^{s}\right)
^{k}\left\vert p_{\mathsf{m}}\right\vert ^{s}\right) ^{\frac{t}{s}}.  \notag
\end{eqnarray}%
Furthermore we have 
\begin{equation}
\dsum\limits_{\mathsf{m}\in \mathbb{Z}}\frac{1}{\left( 1+\left\vert \mathsf{m%
}\right\vert ^{s}\right) ^{\frac{kr}{s}}}<+\infty  \label{convergence}
\end{equation}%
since $kr>r\frac{s-t}{st}=1$ by virtue of the hypothesis and (\ref{duality}%
), and therefore%
\begin{equation*}
\dsum\limits_{\mathsf{m}\in \mathbb{Z}}\hat{w}_{\mathsf{m}}\left\vert p_{%
\mathsf{m}}\right\vert ^{t}\leq c_{k,s,t}\left( \dsum\limits_{\mathsf{m}\in 
\mathbb{Z}}w_{\mathsf{m}}\left( 1+\left\vert \mathsf{m}\right\vert
^{s}\right) ^{k}\left\vert p_{\mathsf{m}}\right\vert ^{s}\right) ^{\frac{t}{s%
}}
\end{equation*}%
for some constant $c_{k,s,t}>0$ depending solely on $k$, $s$ and $t$.
Consequently, changing the value of $c_{k,s,t}$ as necessary we have%
\begin{equation}
\left\Vert \mathsf{p}\right\Vert _{t,\hat{w}}\leq c_{k,s,t}\left\Vert 
\mathsf{p}\right\Vert _{k,s,w}  \label{continuousembedding}
\end{equation}%
for every $\mathsf{p}\in h_{\mathbb{C},w}^{k,s}$, which is the desired
continuity.

As for compactness, if $\mathcal{K}$ is a bounded set in $h_{\mathbb{C}%
,w}^{k,s}$ and if $\kappa >0$ is the radius of a ball centered at the origin
of $h_{\mathbb{C},w}^{k,s}$ and containing $\mathcal{K}$, then for each $%
\varepsilon >0$ there exists $\mathsf{m}^{\ast }\in \mathbb{Z}^{+}$ such that%
\begin{equation*}
\dsum\limits_{\left\vert \mathsf{m}\right\vert \geq \mathsf{m}^{\ast }}\frac{%
1}{\left( 1+\left\vert \mathsf{m}\right\vert ^{s}\right) ^{\frac{kr}{s}}}%
\leq \left( \frac{\varepsilon c_{1}^{\frac{1}{t}}}{2\kappa }\right) ^{r}
\end{equation*}%
because of (\ref{convergence}), where $c_{1}>0$ is chosen to be the constant
that appears in (\ref{weightsbis}) or on the right-hand side of the last
inequality in (\ref{estimates}). Consequently, from an estimate similar to (%
\ref{estimates}) we obtain for every $\mathsf{p}\in \mathcal{K}$ the
inequalities%
\begin{eqnarray*}
&&\dsum\limits_{\left\vert \mathsf{m}\right\vert \geq \mathsf{m}^{\ast }}%
\hat{w}_{\mathsf{m}}\left\vert p_{\mathsf{m}}\right\vert ^{t} \\
&\leq &c_{1}^{-1}\left( \dsum\limits_{\left\vert \mathsf{m}\right\vert \geq 
\mathsf{m}^{\ast }}\frac{1}{\left( 1+\left\vert \mathsf{m}\right\vert
^{s}\right) ^{\frac{kr}{s}}}\right) ^{\frac{t}{r}}\left( \dsum\limits_{%
\mathsf{m\in }\mathbb{Z}}w_{\mathsf{m}}\left( 1+\left\vert \mathsf{m}%
\right\vert ^{s}\right) ^{k}\left\vert p_{\mathsf{m}}\right\vert ^{s}\right)
^{\frac{t}{s}} \\
&\leq &c_{1}^{-1}\left( \frac{\varepsilon c_{1}^{\frac{1}{t}}}{2\kappa }%
\right) ^{t}\kappa ^{t}=\left( \frac{\varepsilon }{2}\right) ^{t},
\end{eqnarray*}%
that is, 
\begin{equation*}
\left( \dsum\limits_{\left\vert \mathsf{m}\right\vert \geq \mathsf{m}^{\ast
}}\hat{w}_{\mathsf{m}}\left\vert p_{\mathsf{m}}\right\vert ^{t}\right) ^{%
\frac{1}{t}}\leq \frac{\varepsilon }{2}.
\end{equation*}%
Compactness of the first embedding in (\ref{embeddings}) then follows from
an argument entirely similar to that set forth in the proof of Theorem 1,
based on the existence of the Schauder basis $\left( \mathsf{e}_{\mathsf{m}%
}\right) _{\mathsf{m\in }\mathbb{Z}}$. The second embedding is an easy
consequence of the second inequality in (\ref{weightsbis}), for we have
successively%
\begin{eqnarray*}
&&\dsum\limits_{\mathsf{m\in }\mathbb{Z}}w_{\mathsf{m}}\left\vert p_{\mathsf{%
m}}\right\vert ^{s} \\
&=&\dsum\limits_{\mathsf{m\in }\mathbb{Z}}\left( w_{\mathsf{m}}^{\frac{t}{s}%
}\left\vert p_{\mathsf{m}}\right\vert ^{t}\right) ^{\frac{s}{t}}\leq \left(
\dsum\limits_{\mathsf{m\in }\mathbb{Z}}w_{\mathsf{m}}^{\frac{t}{s}%
}\left\vert p_{\mathsf{m}}\right\vert ^{t}\right) ^{\frac{s}{t}} \\
&\leq &\left( c_{2}\dsum\limits_{\mathsf{m\in }\mathbb{Z}}\hat{w}_{\mathsf{m}%
}\left\vert p_{\mathsf{m}}\right\vert ^{t}\right) ^{\frac{s}{t}%
}=c_{s,t}\left( \dsum\limits_{\mathsf{m\in }\mathbb{Z}}\hat{w}_{\mathsf{m}%
}\left\vert p_{\mathsf{m}}\right\vert ^{t}\right) ^{\frac{s}{t}}
\end{eqnarray*}%
with an obvious choice for $c_{s,t}>0$. \ \ $\blacksquare $

\bigskip

Let us now consider the particular case of Theorem 2 we alluded to above. We
first define%
\begin{equation*}
\hat{h}_{\mathbb{C},w}^{k,s}:=\left\{ \mathsf{p}\in h_{\mathbb{C},w}^{k,s}:%
\text{ }p_{\mathsf{m}}=0\text{ \ for all }\mathsf{m}\in \mathbb{Z}^{-}\text{ 
}\right\}
\end{equation*}%
and similarly $\hat{l}_{\mathbb{C},\hat{w}}^{t}$ and $\hat{l}_{\mathbb{C}%
,w}^{s}$, which are closed subspaces of $h_{\mathbb{C},w}^{k,s}$, $l_{%
\mathbb{C},\hat{w}}^{t}$ and $l_{\mathbb{C},w}^{s}$, respectively. Then
embeddings (\ref{embeddings}) still hold when\textit{\ }$s>t$ and $k\in
\left( \frac{s-t}{st},+\infty \right) $, that is\ 
\begin{equation*}
\hat{h}_{\mathbb{C},w}^{k,s}\hookrightarrow \hat{l}_{\mathbb{C},\hat{w}%
}^{t}\rightarrow \hat{l}_{\mathbb{C},w}^{s}\text{.}
\end{equation*}%
With $\beta >0$, let $w_{\beta }:=(w_{\beta ,\mathsf{m}})$ be the sequence
of Gibbs related weights given by%
\begin{equation*}
w_{\beta ,\mathsf{m}}=\exp \left[ \beta \mathsf{m}\right]
\end{equation*}%
for every $\mathsf{m}\in \mathbb{N}$ and let $\hat{w}_{\beta }:=w_{\frac{%
\beta }{2}}$. It is then plain that condition (\ref{weightsbis}) holds if we
choose $s=2$, $t=1$ and $k=1$, so that we have%
\begin{equation}
\hat{h}_{\mathbb{C},w_{\beta }}^{1,2}\hookrightarrow \hat{l}_{\mathbb{C},w_{%
\frac{\beta }{2}}}^{1}\rightarrow \hat{l}_{\mathbb{C},w_{\beta }}^{2}
\label{embeddingsbis}
\end{equation}%
where the first embedding is compact and the second one continuous. This
choice of weights and (\ref{embeddingsbis}) are exactly what was used in 
\cite{vuillermot} to prove the self-adjointness and the compactness of the
resolvent of the infinitesimal generator of a class of master equations
describing the interaction of a one-dimensional quantum harmonic oscillator
with a thermal bath at inverse temperature $\beta $. We refer the reader to 
\cite{vuillermot} for details.

\bigskip

\textsc{Remark.} Theorem 3.1 in \cite{bisgard} is also a particular case of
Theorem 2 corresponding to the same values of $s$, $t$ and $k$ as above, but
with $w_{\mathsf{m}}=\hat{w}_{\mathsf{m}}=1$ for every $\mathsf{m}\in 
\mathbb{Z}.$

\bigskip

We complete this section by proving a result of Pitt's type for the spaces $%
h_{\mathbb{C},w}^{k,s}$. Recall that if $C$ denotes any linear bounded
operator from $l_{\mathbb{C}}^{s}$ into $l_{\mathbb{C}}^{t}$ with $s>t\geq 1$%
, where $l_{\mathbb{C}}^{s}$ and $l_{\mathbb{C}}^{t}$ carry the norm (\ref%
{normbis}) with $w_{\mathsf{m}}=1$ for every $\mathsf{m}$, then $C$ is
necessarily compact (see the aforementioned references on this theme).
Remembering that embedding (\ref{continuousembedingbis}) holds we then have:

\bigskip

\textbf{Theorem 3.}\textit{\ Let us assume that }$s>t\geq 1$. \textit{Then
every linear bounded operator }$T:h_{\mathbb{C},w}^{k,s}\rightarrow h_{%
\mathbb{C},w}^{k,t}$\textit{\ is compact.}

\bigskip

\textbf{Proof. }Let us define $J_{k,s,w}:h_{\mathbb{C},w}^{k,s}\rightarrow
l_{\mathbb{C}}^{s}$ by%
\begin{equation*}
\left( J_{k,s,w}\mathsf{p}\right) _{\mathsf{m}}:=w_{\mathsf{m}}^{\frac{1}{s}%
}\left( 1+\left\vert \mathsf{m}\right\vert ^{s}\right) ^{\frac{k}{s}}p_{%
\mathsf{m}}
\end{equation*}%
for every $\mathsf{m}\in \mathbb{Z}$. It is then plain that $J_{k,s,w}$ is
an isometric isomorphism with $J_{k,s,w}\left( h_{\mathbb{C},w}^{k,s}\right)
=$ $l_{\mathbb{C}}^{s}$ and $J_{k,s,w}^{-1}\left( l_{\mathbb{C}}^{s}\right)
=h_{\mathbb{C},w}^{k,s}$, so that%
\begin{equation*}
C:=J_{k,t,w}TJ_{k,s,w}^{-1}
\end{equation*}%
is a linear bounded operator from $l_{\mathbb{C}}^{s}$ into $l_{\mathbb{C}%
}^{t}$, hence compact as a consequence of Pitt's theorem. Therefore%
\begin{equation*}
T=J_{k,t,w}^{-1}CJ_{k,s,w}
\end{equation*}%
is also compact as the composition of $C$ with linear bounded operators. \ \ 
$\blacksquare $

\bigskip

\textsc{Remark.} In view of many potential applications of Orlicz space
theory (see, e.g., \cite{raoren}), an interesting open problem is that of
extending the results of this article to appropriate scales of weighted
sequence spaces of Orlicz-Sobolev type.

\bigskip

\textbf{Statements and declarations}

\bigskip

\textsc{Declaration of competing interest}: The author declares that he has
no known competing financial interests or personal relationships that could
have appeared to influence the work reported in this paper.

\bigskip

\textsc{Conflict of interest:} The author states that there is no conflict
of interest regarding the content of this paper.

\bigskip

\textsc{Data availability statement: }All the data supporting the results
stated in this article are available in the bibliographical references
listed below.

\bigskip

\textsc{Funding.} The author would like to thank the Funda\c{c}\~{a}o para a
Ci\^{e}ncia e a Tecnologia (FCT) of the Portugu\^{e}s Government for its
financial support under grant UIDB/04561/2020.

\end{document}